\documentclass[11pt,oneside]{article}
\usepackage[utf8]{inputenc}
\usepackage[english]{babel}
\usepackage{amsmath}
\usepackage{amsfonts}
\usepackage{amssymb}
\usepackage{amsthm}
\usepackage{comment}
\usepackage{tikz-cd} 
\usepackage{hyperref}
\usepackage{pgfplots}
\usetikzlibrary{calc}
\usepackage[left=3cm,right=3cm,top=3cm,bottom=3cm]{geometry}
\usepackage[title]{appendix}
\usepackage[]{algorithm2e}
\usepackage{comment}
\usepackage[normalem]{ulem}
\useunder{\uline}{\ul}{}
\theoremstyle{plain}
\newtheorem{teo}{}[section]
\newtheorem{prop}[teo]{Proposition}
\newtheorem{cor}[teo]{Corollary}

\newtheorem{lem}[teo]{Lemma}
\newtheorem{thm}[teo]{Theorem}
\theoremstyle{definition}
\newtheorem{ex}[teo]{Example}
\newtheorem{rem}[teo]{Remark}
\newtheorem{df}[teo]{Definition}
\usepackage{graphicx} 

\newcommand\blfootnote[1]{%
  \begingroup
  \renewcommand\thefootnote{}\footnote{#1}%
  \addtocounter{footnote}{-1}%
  \endgroup
}

\title{Compactness-like properties in Alexandroff paratopological groups}
\author{Tayomara Borsich $^{a}$ and Pedro J. Chocano$^{a}$}
\date{}

\begin{document}

\maketitle
\begin{center} {\small $^{a}$ Departamento de Matemática Aplicada, Universidad Rey Juan Carlos. E-mail addresses: tayomara.gonzalez@urjc.es (T. Borsich), pedro.chocano (P.J. Chocano)}
\end{center}

\begin{abstract}
We study compactness-like and covering properties in non-compact Alexandroff paratopological groups. Several properties are characterized in terms of the order structure associated with the Alexandroff topology and the number of connected components. We establish equivalences among $2$-pseudocompactness, countable $o$-compactness, Baire, and an intersection property for increasing chains. We also obtain a common characterization of $\sigma$-local compactness, $\sigma$-compactness, the Lindelöf property, and left $\omega$-precompactness in terms of the existence of a countable cofinal subset. Moreover, we characterize pseudocompactness, $H$-compactness, star-compactness, the countable chain condition, and the weakly Lindelöf property by means of the number of connected components. Finally, we determine further implications among the properties considered and provide examples and counterexamples.
\end{abstract}

\blfootnote{This research is partially supported by Grant  PID2024-156663NB-I00 from Ministerio de Ciencia, Innovación y Universidades (Spain) and  2025/SOLCON-159637 from Rey Juan Carlos University}
\section{Introduction}

In recent years, Alexandroff spaces have attracted increasing attention in various areas of mathematics due to their combinatorial nature and their usefulness in addressing a wide range of problems. For instance, in \cite{barmak2024conley}, finite topological spaces, which constitute a particular class of Alexandroff spaces, are used to develop a suitable theory of the Conley index and apply it to classical problems in dynamical systems. Another interesting example arises in shape theory, a generalization of homotopy theory particularly well suited to the study of compact metric spaces with poor local properties. In \cite{chocano2025shape}, the authors show that finite topological spaces can be used to describe the shape category of compact metric spaces. These are only two representative examples of the numerous applications of Alexandroff spaces that have appeared in the literature; see, for instance, \cite{barmak2007minimal,barmak2020Lefschetz,chocano2025realizingisometries,chocano2026samples,chocanoprieto2026euler}.

In this paper, we focus on groups admitting an Alexandroff topology for which the multiplication map is continuous, that is, on Alexandroff paratopological groups. Our interest in these structures originated in \cite{chocanoPrieto2026negativedimension}, where a group endowed with a natural Alexandroff topology arises. It was proved in \cite{borsich2026paratopological,borsich2027continuous} that every Alexandroff topological group is trivial. Nevertheless, allowing the inversion map to be discontinuous gives rise to nontrivial Alexandroff paratopological groups. However, every compact Alexandroff paratopological group is trivial. Thus, imposing compactness makes the interaction between the group structure and the Alexandroff topology excessively restrictive, leading only to trivial examples. For this reason, in \cite{borsich2026paratopological}, the authors restricted their attention to non-compact Alexandroff paratopological groups. In that paper, they establish the existence of such groups by constructing several examples and investigate some of their basic properties.

Since compactness is a highly restrictive condition in the setting of Alexandroff paratopological groups, in this paper we undertake a systematic study of compactness-like properties and other related notions in this class of groups. We characterize several of these properties and determine the implications and equivalences among them, providing counterexamples whenever the corresponding implications fail. Some of the definitions considered in this manuscript are inspired by \cite{BanakhRavsky2020}, where several properties of general paratopological groups and the relations among them are investigated. To help the reader navigate the large number of notions and relations considered throughout the paper, we summarize our main results in this introduction by means of implication diagrams. Although some of the notions appearing in these diagrams may not yet be familiar to the reader, all the relevant definitions are recalled in the corresponding sections.

Section \ref{sec_cofinal} is organized as follows. In the first part, we analyze and characterize $2$-pseudocompactness, countably $o$-compactness, and the Baire property. In particular, we obtain simple characterizations of these properties in terms of intersections of minimal closed sets associated with increasing chains. More precisely, Theorem \ref{thm_car_2_pseudocompact} establishes the following equivalences for non-compact Alexandroff paratopological groups with finitely many connected components:
\begin{figure}[ht]
\centering 
\[ \begin{tikzcd}[column sep=large,row sep=large] 2\text{-pseudocompactness}  \arrow[d, Leftrightarrow]  \arrow[r, Leftrightarrow] &  \text{countable }o\text{-compactness}  \arrow[d, Leftrightarrow] \\   \text{Every sequence } x_n\leq x_{n+1} \text{ satisifes }\bigcap_{n\in \mathbb{N}} F_{x_n}\neq \emptyset \arrow[r, Leftrightarrow]   & \text{Baire} \end{tikzcd} \]
\caption{Equivalence between $2$-pseudocompactness, countable $o$-compactness and the Baire property.} \label{fig:2pseudo}
\end{figure}

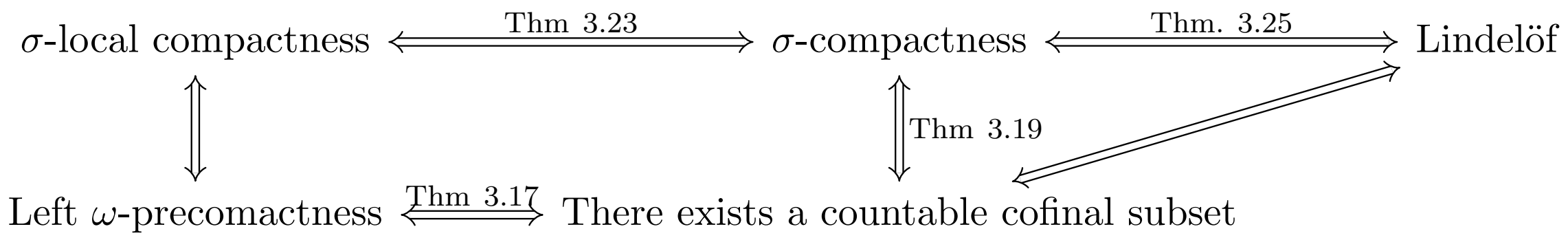
\begin{figure}[ht]
\centering 
\[ \begin{tikzcd}[column sep=large,row sep=large] \sigma\text{-local compactness}  \arrow[r, Leftrightarrow,"\text{Thm}~\ref{thm_sigma_local}"] \arrow[d, Leftrightarrow]   &  \sigma\text{-compactness}\arrow[r, Leftrightarrow,"\text{Thm}.~\ref{thm_sigma_compacto_equivalente_lindelof}"] \arrow[d, Leftrightarrow,"\text{Thm}~\ref{thm_caracterizacion_sigma_compacto}"]   & \text{Lindelöf} \arrow[dl, Leftrightarrow] \\ \text{Left }\omega\text{-precomactness} \arrow[r, Leftrightarrow,"\text{Thm}~\ref{char_left-omega}"]   & \text{There exists a countable cofinal subset} & \end{tikzcd} \]
\caption{Equivalence between  $\sigma$-compactness, $\sigma$-local compactness, Lindelöf and left $\omega$-precompactness.} \label{fig:sigma_compactness}
\end{figure}

In the second part, we establish several equivalences between classical covering properties from general topology, such as the Lindelöf property, and notions specific to paratopological groups, such as left $\omega$-precompactness. In particular, we obtain a simple common characterization of all the properties considered in this part in terms of the existence of a countable cofinal subset. These results are summarized in Figure \ref{fig:sigma_compactness}. The unlabeled equivalences follow by combining the indicated results.

The final part of this section examines the relationship between the two classes described by the diagrams of Figure \ref{fig:2pseudo} and Figure \ref{fig:sigma_compactness}. More precisely, we prove a non-compact Alexandroff paratopological group with finitely many connected components cannot satisfy simultaneously the equivalent properties in the first diagram (see Figure \ref{fig:2pseudo}) and those in the second diagram (see Figure \ref{fig:sigma_compactness}). In particular, $2$-pseudocompactness and $\sigma$-compactness are mutually exclusive in this setting. However, these two classes are not exhaustive, as we provide an example of a connected non-compact Alexandroff paratopological group that belongs to neither of them.

Section \ref{sec_numero_componentes_conexas} is devoted to the classification of non-compact Alexandroff paratopological groups in terms of their number of connected components. In particular, we characterize several classical properties by means of this cardinal invariant. As usual, we denote by $\pi_0(X)$ the set of connected components of a topological space $X$, so that $|\pi_0(X)|$ denotes the number of its connected components. In the first part, we study properties such as pseudocompactness, $H$-compactness, and star-compactness. The results obtained are summarized in the diagram of Figure \ref{fig:pseudocompactness}, where the unlabeled equivalences follow by combining the indicated results:
\begin{figure}[ht]
\centering 
\[ \begin{tikzcd}[column sep=large,row sep=large] \text{Pseudocompactness} \arrow[r, Leftrightarrow] \arrow[dr, Leftrightarrow,"\text{Thm}~\ref{car_pseudocompact}"]   &  H\text{-compactness} \arrow[r, Leftrightarrow] \arrow[d, Leftrightarrow,"\text{Thm}.~\ref{thm_car_h-compact}"]   & \text{Star-compactness} \arrow[dl, Leftrightarrow,"\text{Cor}~\ref{char_star_compact}"] \\  & \lvert \pi_0(X)\rvert <\infty & \end{tikzcd} \]
\caption{Equivalence between pseudocompactness, $H$-compactness and star-compactness.} \label{fig:pseudocompactness}
\end{figure}
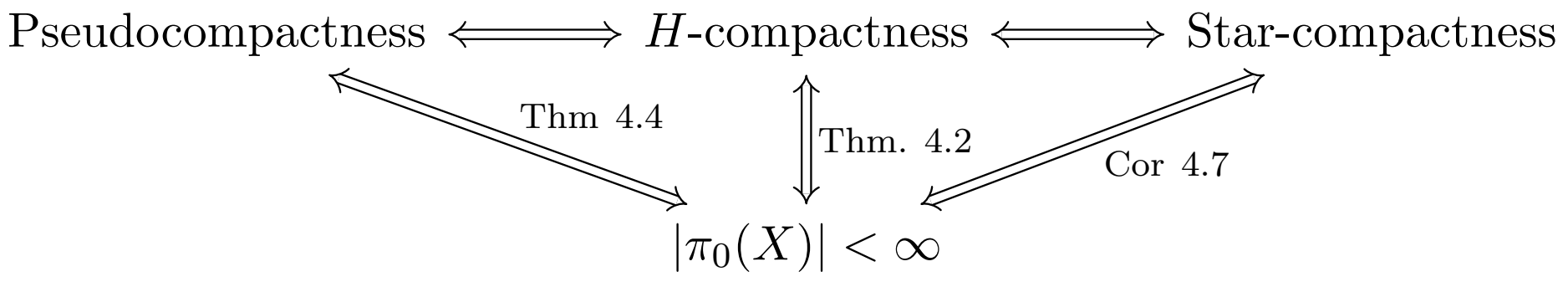

\begin{figure}[ht]
\centering 
\[ \begin{tikzcd}[column sep=large,row sep=large] \text{Countable chain condition}  \arrow[rr, Leftrightarrow] \arrow[dr, Leftrightarrow,"\text{Thm}~\ref{chain_condition}"]  & & \text{Weakly Lindelöf}\arrow[dl, Leftrightarrow,"\text{Thm}~\ref{weakly_lindelof}"]    \\  & \lvert \pi_0(X)\rvert \leq \aleph_0 & \end{tikzcd} \]
\caption{Equivalence between countable chain condition and weakly Lindelöf.} \label{fig:weakly_lindelof}
\end{figure}
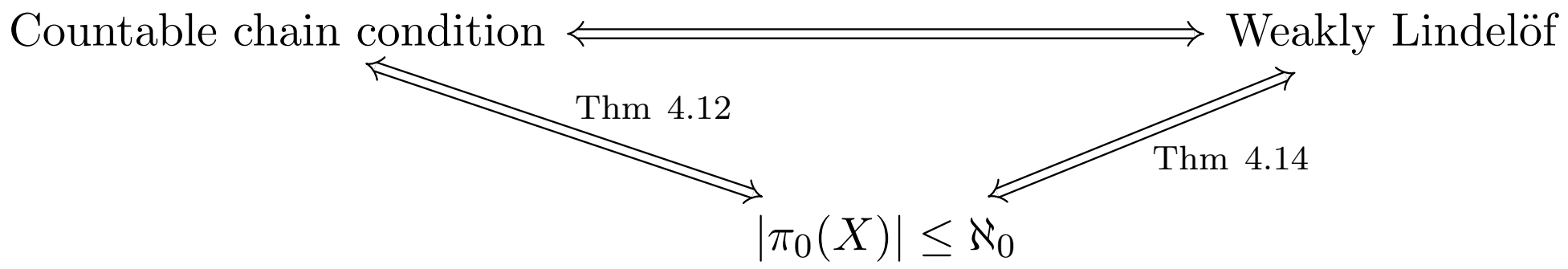

In the second part, we study the countable chain condition and the weakly Lindelöf property. In particular, we prove that both properties are equivalent to having countably many connected components. These equivalences are summarized in the diagram of Figure \ref{fig:weakly_lindelof}.

Note that, under our hypotheses, pseudocompactness, and hence any of the equivalent properties in the diagram of Figure \ref{fig:pseudocompactness}, implies the weakly Lindelöf property, and hence all the equivalent properties in the diagram of Figure \ref{fig:weakly_lindelof}.

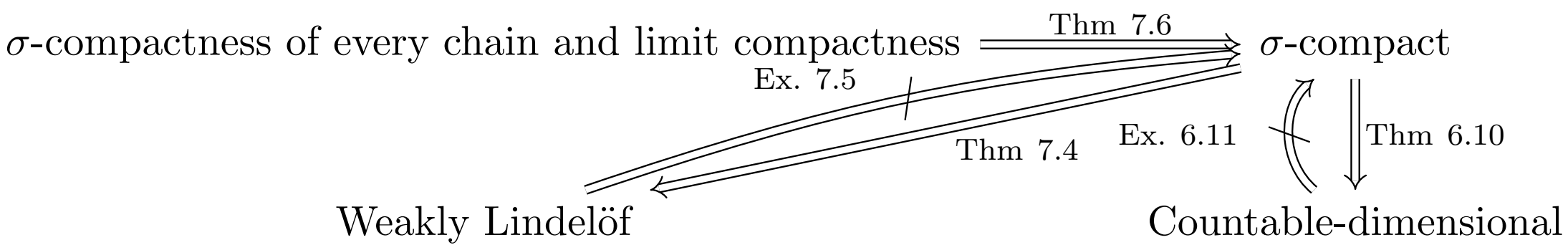
\begin{figure}[ht]
\centering 
\[ \begin{tikzcd}[column sep=large,row sep=large] \sigma\text{-compactness of every chain}\text{ and limit compactness} \arrow[r, Rightarrow,"\text{Thm}~\ref{limit_compact_sigma_compact_chain}"] &\sigma\text{-compact} \arrow[dl, Rightarrow,"\text{Thm}~\ref{sigmacompactoimplicaweaklylindelof}"] \arrow[d, Rightarrow, "\text{Thm}~\ref{sigma_compacto_implica_countable}"]    \\
\text{Weakly Lindelöf} \arrow[ur, Rightarrow,"/" marking, bend left=7, "\text{Ex.}~\ref{ex_weakly_lindelof_no_implica_sigma_compact} \quad"]& \text{Countable-dimensional} \arrow[from=2-2,to=1-2, Rightarrow,"/" marking, bend left=50, "\text{Ex.}~\ref{ex_ZxR} \quad"] \end{tikzcd} \]
\caption{Implications involving $\sigma$-compactness.}
 \label{fig:implicaciones_sigma_compact}
\end{figure}

We summarize in the diagram of Figure \ref{fig:implicaciones_sigma_compact} some additional relations established in Section \ref{sec_covering_dimension} and Section \ref{sec_consequencias}. Particularly, Section \ref{sec_covering_dimension} focuses on the notion of covering dimension within our framework.  For simplicity, we choose one representative from each of the equivalence classes described above.

In Section~\ref{sec_fallan}, we discuss several classical properties that cannot hold in non-compact Alexandroff paratopological groups. More precisely, we show that such groups cannot be $\omega$-bounded, totally countably compact, strongly locally compact, sequentially feebly compact, precompact or metacompact.

We conclude this introduction with a brief description of the organization of the paper. Section \ref{sec_preliminaries} is devoted to preliminaries. We recall the basic notions and results concerning Alexandroff spaces and paratopological groups, together with some elementary properties of Alexandroff paratopological groups established in \cite{borsich2026paratopological} that will be used throughout the paper. The subsequent sections contain the results summarized in the diagrams above.

\section{Preliminaries}\label{sec_preliminaries}

We recall the main definitions from the theory of Alexandroff topological spaces and also state notation for future sections. For a detailed introduction to this theory we recommend \cite{barmak2011algebraic,may1966finite}. 

\begin{df}
    An \emph{Alexandroff space} is a topological space in which arbitrary intersections of open sets are open.
\end{df}

Given an Alexandroff space $X$ and $x\in X$, define $U_x$ (resp. $F_x$) as the intersection of every open (closed) set containing the point $x$. Note that $U_x$ (resp. $F_x$) is the minimal open (resp. closed) set containing $x$.  Consider the relation $\leq$ on $X$ defined by $x\leq y$ if and only if $U_x\subseteq U_y$. It is straightforward to verify that $X$ with this relation is a preordered set. Moreover, if $X$ is a $T_0$ space, then $(X,\leq )$ is a partially ordered set (or poset for short).

Given a preordered set $(X,\leq)$, the set of upper (resp. lower) sets of $X$ forms a basis for an Alexandroff topology on $X$. Particularly, if $(X,\leq)$ is a poset, then the topological space $X$ is $T_0$.

These two construction are mutually inverse and we also have that any map $f:X\rightarrow Y$ between two Alexandroff spaces is continuous if and only if $f:(X,\leq_X)\rightarrow (Y,\leq_Y)$ is order preserving. From this, 

\begin{thm}
    The category of preordered (resp. partially ordered) sets and the category of Alexandroff topological (resp. $T_0$)  spaces are isomorphic.
\end{thm}

From now on, we will not distinguish between Alexandroff spaces or preordered sets and we will combine notation from both theories. For instance, $U_x=\{y\in X : y\leq x\}$ and $F_x=\{y\in X : y\geq x\}$. Particularly, upper sets correspond to closed sets, while lower sets correspond to open sets.





Given an Alexandroff space \(X\), the \emph{width} of \(X\), denoted \(\mathrm{width}(X)\), is defined as the cardinality of a maximal antichain in \(X\). We conclude this part of the section by proving a simple result concerning compactness.
\begin{thm}\label{thm_caracterizacion_compactos}
    Let $X$ be an Alexandroff space. Then $X$ is compact if and only if the set of maximal points of $X$ is finite.
\end{thm}
\begin{proof}
   Let $X$ be a compact space. Suppose, towards a contradiction, that the set of maximal points of $X$ is infinite. Let $M$ denote this set. Then \[ \mathcal U=\{U_x:x\in M\} \] is an open cover of $X$. Since distinct maximal points do not belong to each other's minimal neighborhoods, no finite subfamily of $\mathcal U$ covers $X$. This contradicts the compactness of $X$. Conversely, suppose that the set of maximal points of $X$ is finite, and let $\mathcal U$ be an open cover of $X$. For each maximal point $x$, choose $V_x\in\mathcal U$ such that $x\in V_x$. Since $x$ is maximal, we have $U_x\subseteq V_x$. Let $M$ be the set of maximal points of $X$. Then \[ \mathcal V=\{V_x:x\in M\} \] is a finite family. Moreover, \[ X=\bigcup_{x\in M}U_x\subseteq\bigcup_{x\in M}V_x. \] Therefore, $\mathcal V$ is a finite subcover of $\mathcal U$, and hence $X$ is compact.
\end{proof}
\begin{cor}\label{cor_U_x_compactness}
    Let $X$ be an Alexandroff space and $x\in X$. Then $U_x$ is compact. Particularly, each Alexandroff space is locally compact.
\end{cor}

\begin{df}
    Let $G$ be a group endowed with a topology. We say that $G$ is a \emph{topological group} if:
    \begin{enumerate}
        \item the multiplication map $m\colon G\times G\to G$ is continuous, and
        \item the inversion map $\mathrm{Inv}\colon G\to G$ is continuous.
    \end{enumerate}
    We say that $G$ is a \emph{paratopological group} if only the first condition is required.
\end{df}

In \cite{borsich2026paratopological}, the authors study topological groups endowed with Alexandroff topologies. One of the main results of that paper shows that, under the assumption of being an Alexandroff space, only trivial topological groups can arise. In contrast, there are no trivial examples among non-compact Alexandroff paratopological groups. For this reason, we restrict our study to the non-compact case. We conclude this section by recalling some basic terminology and results from \cite{borsich2026paratopological} that will be used throughout the remainder of the paper.

\begin{prop}\label{propiedadesbasicas} Let $X$ be a non-compact Alexandroff paratopological group. Then the following properties hold: \begin{enumerate} 
\item $X$ does not have maximal points. 
\item If $U$ is open (resp.\ closed), then $U^{-1}$ is closed (resp.\ open). 
\item For every $x\in X$, $U_x^{-1}=F_{x^{-1}}$.
\item Every continuous map $f\colon X\to\mathbb{R}$, where $\mathbb{R}$ is endowed with its usual topology, is constant on each connected component of $X$.
\end{enumerate} \end{prop}

\begin{df} 
A topological space $X$ is said to be {\it hyperconnected} (resp. {\it ultraconnected}) if the intersection of every pair of nonempty open sets (resp. closed sets) is nonempty.
Equivalently, $X$ has no two disjoint nonempty open subsets (resp. closed subsets).
\end{df}

\begin{thm}\label{every_paratopological_alexandroff_hyperconexosiiiiii}\cite[Theorem 4.13]{borsich2026paratopological}
   Every connected non-compact Alexandroff paratopological group $X$ is hyperconnected and ultraconnected.
\end{thm}
\begin{thm}\label{dense}
 Let $X$ be a non-compact Alexandroff paratopological group, then every open subset of $X$ is dense if and only if $X$ is connected.
\end{thm}
\begin{proof}
We argue by contradiction. Suppose that every nonempty open subset of \(X\) is dense and that \(X\) is not connected. Then \(X\) admits a separation \(X=U\cup V\), where \(U\) and \(V\) are disjoint nonempty open subsets of \(X\). Since \(U\) is dense and \(V\) is nonempty and open, it follows that \(U\cap V\neq\emptyset\), a contradiction.

Suppose $X$ is connected. By Theorem \ref{every_paratopological_alexandroff_hyperconexosiiiiii}, $X$ is hyperconnected and, consequently, every nonempty open set is dense in X.

\end{proof}
\begin{cor}\label{densdisc}
 Let $X$ be a non-compact Alexandroff paratopological group. Then every open set that intersects each connected component of $X$ is dense.
\end{cor}
\begin{proof}
The proof is straightforward from Theorem \ref{dense}.
\end{proof}
The following result collects several characterizations of connectedness  within this context.
\begin{thm}\label{conexion_caracterizacion} Let $X$ be a non-compact Alexandroff paratopological group. Then the following properties are equivalent:
\begin{enumerate}
    \item $X$ is connected.
    \item $X$ is path-connected.
    \item Every open subset of $X$ is dense.
    \item $X$ is hyperconnected.
    \item $X$ is ultraconnected.
    \item $X$ is a directed set.
    \item $X$ is a codirected set.
\end{enumerate}
\end{thm}
Let us recall that a directed (resp. codirected) set $X$ is a preordered set satisying that for every $x,y\in X$ there exists $z\in X$ such that $x,y\leq z$ (resp. $x,y\geq z$).
\begin{rem}
    Given an Alexandroff paratopological group $X$ and $x\in X$, it is simple to deduce that $U_x$ is connected, path-connected, hyperconnected and ultraconnected. Hence, $X$ is locally  connected, path-connected, hyperconnected and ultraconnected.
\end{rem}

To conclude this section, and for the sake of completeness, we study one of the remaining classical connectedness properties in this setting.

\begin{df}
A topological space $X$ is said to be arcwise connected (or arc-connected) if any two points of $X$ can be joined by an arc. Recall that an arc joining $x,y\in X$ is a continuous map $f\colon [0,1]\to X$ that is a homeomorphism onto its image and satisfies $f(0)=x$ and $f(1)=y$.
\end{df}

\begin{prop}
    Let $X$ be a non-compact Alexandroff paratopological group. Then $X$ is not arcwise connected.
\end{prop}
\begin{proof}
   As an immediate consequence of \cite{borsich2027continuous}, every continuous map $f\colon X\to Y$ is constant whenever $Y$ is a $T_1$-space. From this, we derive the desired result.
\end{proof}

\section{Classification by cofinal subsets and increasing sequences}\label{sec_cofinal}
We devote this section to a first classification of connected non-compact Alexandroff paratopological groups in terms of compact-like properties. In particular, we study the notions of $2$-pseudocompactness and $\sigma$-compactness and show that they are mutually exclusive within this setting. Moreover, we provide several equivalent conditions for each of these properties.

\subsection{$2$-pseudocompactness, countably $o$-compactness and the Baire property}
In this subsection, we characterize $o$-compactness, 2-pseudocmpactness and the Baire property. Particularly, we demonstrate that they are equivalent notions within the theory of non-compact Alexandroff paratopological groups with finitely many connected components. We start by recalling these notions.
\begin{df}
A paratopological group $X$ is called \emph{$2$-pseudocompact} if, for every decreasing sequence $\{U_n\}_{n\in\mathbb{N}}$ of nonempty open subsets of $X$, one has \[ \bigcap_{n\in\mathbb{N}} \overline{U_n^{-1}} \neq \emptyset. \]
\end{df}

    \begin{df}
    A topological space $X$ is {\it countably $o$-compact} if each decreasing sequence $\{U_n\}_{n\in\mathbb{N}}$ of nonempty open sets in $X$ has nonempty intersection
$\bigcap_{n\in\mathbb{N}} U_n$. 
\end{df}

\begin{df}
    A topological space $X$ is {\it Baire}  if the
intersection of every decreasing sequence of dense open sets is a dense set. 
\end{df}
Before turning to the setting of non-compact Alexandroff paratopological groups, we recall and prove some properties of the above definitions for paratopological groups with finitely many connected components.
\begin{thm}\label{thm_2-pseudo_finitos}
    Let $X$ be a paratopological group with finitely many connected components. Then
    \begin{enumerate}
    \item $X$ is $2$-pseudocompact if and only if each connected component is $2$-pseudocompact.
    \item $X$ is countably $o$-compact if and only if each connected component is countably $o$-compact.
    \item $X$ is Baire if and only if each connected component is Baire.
    \end{enumerate}
    \end{thm}

    \begin{proof}
    1.  If $X$ is $2$-pseudocompact, then the result holds trivially. Suppose that every connected component of $X$ is $2$-pseudocompact. Let $\{U_n\}_{n\in\mathbb N}$ be a decreasing sequence of nonempty open subsets of $X$. For every $n\in\mathbb N$, let \[ \mathcal C_n = \{C:C \text{ is a connected component of }X \text{ and }U_n\cap C\neq\emptyset\}. \] Each $\mathcal C_n$ is finite and nonempty. Moreover, since $U_{n+1}\subseteq U_n$, we have $\mathcal C_{n+1}\subseteq\mathcal C_n$. Thus, the decreasing sequence $\{\mathcal C_n\}_{n\in\mathbb N}$ eventually stabilizes, and hence $\bigcap_{n\in\mathbb N}\mathcal C_n\neq\emptyset$. Choose a connected component $C\in\bigcap_{n\in\mathbb N}\mathcal C_n$.  Then $V_n=U_n\cap C$ is a nonempty open subset of $C$ for every $n\in\mathbb N$, and the sequence $\{V_n\}_{n\in\mathbb N}$ is decreasing. Since $C$ is $2$-pseudocompact, \[ \bigcap_{n\in\mathbb N}\overline{V_n^{-1}} \neq\emptyset. \] Since $V_n\subseteq U_n$, we have $\overline{V_n^{-1}} \subseteq \overline{U_n^{-1}}$ for every $n\in\mathbb N$. Therefore, \[ \bigcap_{n\in\mathbb N}\overline{U_n^{-1}} \neq\emptyset, \] and hence $X$ is $2$-pseudocompact.

    2. The corresponding result for countable $o$-compactness follows from an analogous argument, so we omit the proof.

    3. Since $X$ has finitely many connected components, they are all open. Thus, $X$ is the topological sum of its connected components. The result now follows from the fact that every open subspace of a Baire space is a Baire space and that a topological sum of Baire spaces is a Baire space.

    \end{proof}

Now, we obtain useful properties that will be used in the subsequent characterizations.
\begin{thm}\label{thm_invierte_orden_inversa}
    Let $X$ be a non-compact Alexandroff paratopological group and let $x,y\in X$. Then $x\leq y$ if and only if $x^{-1}\geq y^{-1}$.
\end{thm}
\begin{proof}
    Suppose that $x\leq y$. Since the multiplication map is continuous, $e\leq x^{-1}y$, where $e$ denotes the identity element, and, therefore, $y^{-1}\leq x^{-1}$. The other implication follows similarly.
    
\end{proof}
\begin{cor}
    Let $X$ be a non-compact Alexandroff paratopological group and $\{x_n\}_{n\in\mathbb{N}}\subseteq X$. Then $x_n\leq x_{n+1}$ for every $n\in \mathbb{N}$ if and only if $x_{n}^{-1}\geq x_{n+1}^{-1}$ for every $n\in \mathbb{N}$.
\end{cor}

\begin{cor}\label{interseccionnovacia}
     Let $X$ be a non-compact Alexandroff paratopological group and $\{x_n\}_{n\in\mathbb{N}}\subseteq X$. 
    Then    $\bigcap_{n\in \mathbb{N}} F_{x_n^{-1}}\neq \emptyset$ if and only if $\bigcap_{n\in \mathbb{N}} U_{x_n}\neq \emptyset$.
\end{cor}
\begin{proof}
    Suppose that $\bigcap_{n\in \mathbb{N}} F_{x_n^{-1}}\neq \emptyset$. If $y\in \bigcap_{n\in \mathbb{N}} F_{x_n^{-1}}$, then $y\geq x_n^{-1}$ for every $n\in \mathbb{N}$. By Theorem \ref{thm_invierte_orden_inversa}, $x_n\geq y^{-1}$ for every $n\in \mathbb{N}$, which gives that $y^{-1}\in \bigcap_{n\in \mathbb{N}} U_{x_n}$. The other implication is similar.
\end{proof}

\begin{lem}\label{lem_sucesion_decreciente_existe_cadena}
    Let $X$ be a connected non-compact Alexandroff paratopological group. If $\{U_n \}_{n\in \mathbb{N}}$ (resp. $\{F_n \}_{n\in \mathbb{N}}$) is a decreasing sequence of nonempty open (resp. closed) sets, then there exists $\{x_n \}_{n\in \mathbb{N}}\subset X$ such that $x_n\in U_n$ and $x_{n}\geq x_{n+1}$ (resp. $x_n\in F_n$ and $x_{n}\leq x_{n+1} $) for every $n\in \mathbb{N}$.
\end{lem}
\begin{proof}
   We construct the sequence $\{x_n\}_{n\in\mathbb{N}}\subseteq X$ recursively. Suppose that $x_n$ has already been constructed. Choose $y_{n+1}\in U_{n+1}$ (resp. $y_{n+1}\in F_{n+1}$ ). By Theorem \ref{conexion_caracterizacion}, $X$ is hyperconnected (resp. ultraconnected), so there exists $x_{n+1}\in X$ such that $x_{n+1}\leq y_{n+1}$ and $x_{n+1}\leq x_n$ (resp.  $x_{n+1}\geq y_{n+1}$ and $x_{n+1}\geq x_n$ ). Since $U_{n+1}$ is open (resp. $F_{n+1}$ is closed), and hence a lower set (resp. an upper set), it follows that $x_{n+1}\in U_{n+1}$ (resp. $x_{n+1}\in F_{n+1}$). 
\end{proof}

\begin{thm}\label{thm_car_2_pseudocompact}
    Let $X$ be a non-compact Alexandroff paratopological group with finitely many connected components. Then the following properties are equivalent:
    \begin{enumerate}
        \item $X$ is $2$-pseudocompact.
        \item $X$ is countably $o$-compact.
        \item $X$ is Baire.
        \item For every sequence $\{x_n \}_{n\in \mathbb{N}}\subseteq X$ such that $x_{n+1}\leq x_{n}$ (resp. $x_{n+1}\geq x_n$), where $n\in \mathbb{N}$, $\bigcap_{n\in \mathbb{N}} F_{x_n^{-1}}\neq \emptyset$ (resp. $\bigcap_{n\in \mathbb{N}} F_{x_n}\neq \emptyset$).
        \item For every sequence $\{x_n \}_{n\in \mathbb{N}}\subseteq X$ such that $x_{n+1}\leq x_{n}$ (resp. $x_{n+1}\geq x_n$), where $n\in \mathbb{N}$, $\bigcap_{n\in \mathbb{N}} U_{x_n}\neq \emptyset$ (resp. $\bigcap_{n\in \mathbb{N}} U_{x_n^{-1}}\neq \emptyset$).
    \end{enumerate}
\end{thm}
\begin{proof}

It suffices to prove the result in the connected case, since the general case follows immediately from Theorem \ref{thm_2-pseudo_finitos}. Henceforth, we assume that $X$ is connected.

The equivalence between $4$ and $5$ follows immediately from Corollary \ref{interseccionnovacia}. By Lemma \ref{lem_sucesion_decreciente_existe_cadena}, $5$ implies $2$. Moreover, the equivalence between $2$ and $3$ follows directly from Theorem \ref{conexion_caracterizacion}.

We prove that 2 implies 1. By hypothesis, each decreasing sequence $\{U_n\}_{n\in\mathbb{N}}$ of nonempty open sets in $X$ has nonempty intersection
$\bigcap_{n\in\mathbb{N}} U_n$. Then there exists $t\in \bigcap_{n\in\mathbb{N}} U_n$ such that $t\in U_n$ for every $n\in \mathbb{N}$. Hence, $t^{-1}\in U^{-1}_n$  for every $n\in \mathbb{N}$. By Proposition \ref{propiedadesbasicas}, $U^{-1}_n$ is a closed set and we have that $\overline{U^{-1}_n}=U^{-1}_n$. Thus, $t^{-1}\in \bigcap_{n\in \mathbb{N}} \overline{U^{-1}_n}$.


    Finally, we demonstrate that 1 is equivalent to 4. Suppose that $X$ is $2$-pseudocompact. Let $\{x_n \}_{n\in \mathbb{N}}\subseteq X$ such that $x_{n+1}\leq x_{n}$ for every $n\in \mathbb{N}$. Therefore, $\{ U_{x_n}\}_{n\in \mathbb{N}}$ is a decreasing sequence of nonempty open subsets of $X$. By Proposition \ref{propiedadesbasicas}, $U_{x_n}^{-1}=F_{x_n^{-1}}$. Since $X$ is $2$-pseudocompact, we deduce that $\bigcap_{n\in \mathbb{N}} F_{x_n^{-1}}\neq \emptyset$.

    Now we prove the other implication. Assume that for every subset $\{x_n \}_{n\in \mathbb{N}}\subseteq X$ such that $x_{n+1}\leq x_{n}$ for every $n\in \mathbb{N}$, one has $\bigcap_{n\in \mathbb{N}} F_{x_n^{-1}}\neq \emptyset$. Let $\{U_{n} \}_{n\in \mathbb{N}}$ be a decreasing sequence of nonempty open sets. By Lemma \ref{lem_sucesion_decreciente_existe_cadena}, there exists a sequence $\{x_n \}_{n\in \mathbb{N}}$ such that $x_n\in U_n$ and $x_{n+1}\le x_n$ for every $n\in \mathbb{N}$. Moreover, $U_{x_n}\subseteq U_n$ for every $n\in \mathbb{N}$. By Proposition \ref{propiedadesbasicas},  $\emptyset \neq \bigcap_{n\in \mathbb{N}} F_{x_n^{-1}}\subseteq \bigcap_{n\in \mathbb{N}}\overline{{U_n^{-1}}}$.

\end{proof}

\begin{rem}
    The assumption that $X$ has finitely many connected components is essential. Indeed, let $Y$ be a connected non-compact $2$-pseudocompact Alexandroff paratopological group, and let $\mathbb{Z}$ be endowed with the discrete topology. Then $Y\times\mathbb{Z}$ is a non-compact Alexandroff paratopological group whose connected components are the sets $Y\times\{m\}$, $m\in\mathbb{Z}$. In particular, it has infinitely many connected components, each of which is $2$-pseudocompact. Nevertheless, the sequence of nonempty open sets $U_n=Y\times[n,\infty)$, $n\in\mathbb{N}$, is decreasing and satisfies \[ \bigcap_{n\in\mathbb{N}}\overline{U_n^{-1}} =\bigcap_{n\in\mathbb{N}}Y\times(-\infty,-n] =\emptyset. \] Thus, $Y\times\mathbb{Z}$ is not $2$-pseudocompact. Similarly, we can use also the same sequence $\{U_n\}_{n\in \mathbb{N}}$ to prove that $Y\times \mathbb{Z}$ is neither countably $o$-compact nor Baire.
\end{rem}

\begin{cor}
    Let $X$ be a non-compact Alexandroff paratopological group with finitely many connected components. If $X$ is not $2$-pseudocompact, then there exists an increasing sequence $\{ F_n\}_{n\in \mathbb{N}}$ (resp. $\{ U_n\}_{n\in \mathbb{N}}$) of nonempty proper closed (resp. open) subsets of $X$ such that $X=\bigcup_{n\in \mathbb{N}} F_n$ (resp. $X=\bigcup_{n\in \mathbb{N}} U_n$).
\end{cor}
\begin{proof}
    By Theorem \ref{thm_car_2_pseudocompact}, there exists a decreasing sequence $\{U_{x_n}\}_{n\in \mathbb{N}}$, where $x_{n}\geq x_{n+1}$ for every $n\in \mathbb{N}$, such that $\bigcap_{n\in \mathbb{N}}U_{x_n}=\emptyset$. By considering the complementary, we obtain \[ \bigcup_{n\in \mathbb{N}} (X\setminus U_{x_n})=X.\]
    Let $F_n$ denote $X\setminus U_{x_n}$. Then $\{F_n\}_{n\in \mathbb{N}}$ is an increasing sequence of nonempty proper closed subsets. The remaining part of the proof follows analogously. 
\end{proof}

We now present two non-trivial examples. The first fails to be $2$-pseudocompact, whereas the second is $2$-pseudocompact.

\begin{ex}\label{ex_z_sigma_compacto}
    Let us consider the group of the integer numbers $(\mathbb{Z},+)$ with its usual order. This is a non-compact Alexandroff paratopological group. Consider $\{ -n\}_{n\in \mathbb{N}}\subset \mathbb{Z}$. It is evident that $\bigcap_{n\in \mathbb{N}} U_{n}=\bigcap_{n\in \mathbb{N}}(-\infty,-n]=\emptyset$. By Theorem \ref{thm_car_2_pseudocompact}, we deduce that $\mathbb{Z}$ is not $2$-pseudocompact.
\end{ex}

\begin{ex}\label{ex_2-pseudocompact}
    Let $\mathbb{Z}^{\mathbb{N}}$ denote the set of all maps $f\colon\mathbb{N}\to\mathbb{Z}$. Define a relation $\leq$ on $\mathbb{Z}^{\mathbb{N}}$ by
    \[ f\leq g \quad\Longleftrightarrow\quad \text{there exists }N\in\mathbb{N}\text{ such that } f(n)\leq g(n)\text{ for every }n\geq N. \] 
    
    It is straightforward to verify that $(\mathbb{Z}^{\mathbb{N}},\leq)$ is a preordered set. Endow $\mathbb{Z}^{\mathbb{N}}$ with pointwise addition, defined by 
    
    \[ (f+g)(n)=f(n)+g(n) \] 
    
    for every $f,g\in\mathbb{Z}^{\mathbb{N}}$ and $n\in\mathbb{N}$. Then $(\mathbb{Z}^{\mathbb{N}},+)$ is an Abelian group, and addition is monotone with respect to the product preorder. More precisely, if $f\leq f'$ and $g\leq g'$, then \[ f+g\leq f'+g'. \]

    Thus, $\mathbb{Z}^\mathbb{N}$ is a non-compact Alexandroff paratopological group that is not $T_0$. We will use the characterization provided in Theorem \ref{thm_car_2_pseudocompact} to prove that $\mathbb{Z}^\mathbb{N}$ is $2$-pseudocompact. Consider a sequence of points $\{f_n\}_{n\in \mathbb{N}}\subset \mathbb{Z}^\mathbb{N}$ such that $f_n\geq f_{n+1}$ for every $n\in \mathbb{N}$. Define $g:\mathbb{N}\rightarrow \mathbb{Z}$ by 
    \[
    g(n)=\min \{f_1(n),\ldots, f_n(n) \}-1.\] 
    For each \(k\in\mathbb{N}\) and every \(n\geq k\), we have \[ g(n)=\min\{f_1(n),\ldots,f_n(n)\}-1<f_k(n). \] Therefore, \(g<f_k\) for every \(k\in\mathbb{N}\). Thus, $g\in U_{f_n}$ for every $n\in \mathbb{N}$ and we conclude that $\mathbb{Z}^\mathbb{N}$ is $2$-pseudocompact.


\end{ex}
\begin{rem}
    Example~\ref{ex_2-pseudocompact} can also be formulated in the setting of  non-compact $T_0$ Alexandroff paratopological groups. Indeed, identify two elements of $\mathbb{Z}^{\mathbb{N}}$ whenever they are eventually equal, or equivalently, define 
    \[ f\sim g \quad\Longleftrightarrow\quad f\leq g \text{ and } g\leq f. \] 
    The addition on $\mathbb{Z}^{\mathbb{N}}$ induces a well-defined group operation on the quotient $\mathbb{Z}^{\mathbb{N}}/{\sim}$, and the induced partial order is compatible with this operation. Moreover, the argument from Example~\ref{ex_2-pseudocompact} shows that every sequence in the quotient admits a strict lower bound which gives that $\mathbb{Z}^{\mathbb{N}}/{\sim}$ is $2$-pseudocompact.
\end{rem}

\subsection{Left $\omega$-precompactness, $\sigma$-compactness, $\sigma$-locally compactness and Lindelöf}

We recall the main definitions of this subsection and then provide characterizations.

\begin{df}\label{def_left_omega_precompact}
    A paratopological group $G$ is {\it left $\omega$-precompact}  if for each neighborhood $U$ of the identity of $G$ there exists a countable subset $C$ of $G$ such that $CU=G$.
\end{df}

 \begin{lem}\label{lem_char_w_narrow}
    Let $X$ be a non-compact Alexandroff paratopological group. Then $X$ is left $\omega$-precompact if and only if there exists a countable set $C\subseteq X$ such that $X=CU_e$, where $e$ denotes the identity element of $X$.
\end{lem}
\begin{proof}
    The proof is immediate from Definition \ref{def_left_omega_precompact} and the fact that every open set $U$ containing the identity element $e$ satisfies that $U_e\subseteq U$.
\end{proof}
Let us recall that given any preordered set $X$, a cofinal subset $C\subseteq X$ is a set satisfying that for every $x\in X$ there exists $c_x\in C $ such that $x\leq c_x$.
\begin{thm}\label{char_left-omega}
    Let $X$ be a non-compact Alexandroff paratopological group. Then $X$ is left $\omega$-precompact if and only if there exists a countable cofinal subset in $X$. 
\end{thm}
\begin{proof}
 Let $e$ denote the identity element of $X$. If $C$ is a countable cofinal subset in $X$, then for every $x\in X$ there exists $c_x\in C$ such that $x\leq c_x$. This yields that $c_x^{-1}x\leq e$, i.e., $c_x^{-1}x\in U_e$. Hence, $x=c_xc_x^{-1}x\in CU_e$ and we conclude that $CU_e=X$. Finally, by Lemma \ref{lem_char_w_narrow}, we obtain the desired result.

 Suppose now that $X$ is left $\omega$-precompact. Then, by Lemma \ref{lem_char_w_narrow}, there exists a countable subset $C$ such that $CU_e=X$. Let us verify that $C$ is a cofinal subset of $X$. If $x\in X$, then there exist $c_x\in C$ and $z\in U_e$ such that $x=c_xz$ or, equivalently, $xc_x^{-1}=z$. Since $z\in U_e$, we have that $xc_x^{-1}\leq e$. By the continuity of the multiplication map, we obtain that $x\leq c_x$. This proves that $C$ is a cofinal subset of $X$.
\end{proof}

\begin{df}
    A topological space $X$ is called {\it $\sigma$-compact} if it is the union of countably many compact sets.
\end{df}

\begin{thm}\label{thm_caracterizacion_sigma_compacto}
    Let $X$ be a non-compact Alexandroff paratopological group. Then $X$ is $\sigma$-compact if and only if there exists a countable cofinal subset in $X$. 
\end{thm}
\begin{proof}
    Suppose that $X$ is $\sigma$-compact. Then $X=\bigcup_{n\in \mathbb{N}}C_n$ where $C_n$ is compact for every $n\in \mathbb{N}$. By Theorem \ref{thm_caracterizacion_compactos}, $C_n$ has a finite number of maximal points. Let $M_n$ denote the set of maximal points of $C_n$. Then $\bigcup_{n\in \mathbb{N}} M_n$ is a countable cofinal subset in $X$.

    If there exists a countable cofinal subset $M$ in $X$, then $X=\bigcup_{x\in M} U_x$. By Corollary \ref{cor_U_x_compactness}, $U_x$ is compact and from this we conclude that $X$ is $\sigma$-compact.

\end{proof}
\begin{cor}\label{cor_numerable_implica_sigma_compacto}
    Let $X$ be a non-compact Alexandroff paratopological group. If $X$ is countable, then $X$ is $\sigma$-compact. 
\end{cor}

\begin{ex}
    Consider the paratopological group $\mathbb{Z}$ from Example \ref{ex_z_sigma_compacto}. By Corollary \ref{cor_numerable_implica_sigma_compacto}, $\mathbb{Z}$ is $\sigma$-compact.
\end{ex}

\begin{df}
    A topological space $X$ is called {\it $\sigma$-locally compact} if it is $\sigma$-compact and locally compact.
\end{df}

\begin{thm}\label{thm_sigma_local}
    Let $X$ be a non-compact Alexandroff paratopological group. Then $X$ is $\sigma$-locally compact if and only if $X$ is $\sigma$-compact.
\end{thm}
\begin{proof}
    It is an immediate consequence of Corollary \ref{cor_U_x_compactness}.
\end{proof}


\begin{df}
    A topological space $X$ is
called {\it Lindelöf} if every open cover has a countable subcover.
\end{df}

\begin{thm}\label{thm_sigma_compacto_equivalente_lindelof}
    Let $X$ be a non-compact Alexandroff paratopological group. Then $X$ is $\sigma$-compact if and only if $X$ is Lindelöf.
\end{thm}
\begin{proof}
    Suppose that $X$ is $\sigma$-compact. By Theorem \ref{thm_caracterizacion_sigma_compacto}, $X$ has a countable cofinal subset $M$. Let $\mathcal{U}$ be an open cover of $X$. For each point $x\in M$, there exists at least one open set $V_x\in \mathcal{U}$ containing $x$. Consider $\mathcal{V}=\{V_x \}_{x\in M}$. This is a subcover that is countable, as desired.
    
    Now, suppose that $X$ is Lindelöf. Then every open cover $\mathcal{U}$ of $X$ contains a countable subcover $\mathcal{V}$. This gives that the open cover $\mathcal{U}=\{U_x \}_{x\in X}$ contains a countable subcover, i.e., $\mathcal{V}=\{U_x \}_{x\in M}$ where $M$ is a countable subset of $X$. By the definition of cover and $U_x$, we may conclude that $M$ is a countable cofinal subset of $X$. Hence, $X$ is $\sigma$-compact by Theorem \ref{thm_caracterizacion_sigma_compacto}.
\end{proof}

By combining Theorems \ref{char_left-omega}, \ref{thm_caracterizacion_sigma_compacto}, \ref{thm_sigma_local} and \ref{thm_sigma_compacto_equivalente_lindelof}, we obtain the following characterization:
\begin{thm}\label{char_sigma_compact}
    Let $X$ be a non-compact Alexandroff paratopological group. Then the following properties are equivalent:
    \begin{enumerate}
        \item $X$ is $\sigma$-compact.
        \item $X$ is $\sigma$-locally compact.
        \item $X$ is Lindelöf.
        \item $X$ is left $\omega$-precompact.
        \item $X$ has a countable cofinal subset.
    \end{enumerate}
\end{thm}

\subsection{Main result}
This subsection contains the main result of the section, which shows that $2$-pseudocompactness and $\sigma$-compactness are mutually exclusive for connected non-compact Alexandroff paratopological groups.

\begin{thm}\label{thm_excluyentes}
    Let $X$ be a connected non-compact Alexandroff paratopological group. Then $X$ cannot be both $2$-pseudocompact and $\sigma$-compact.
\end{thm}

\begin{proof}

     Let us assume that $X$ is $2$-pseudocompact. By Theorem \ref{thm_car_2_pseudocompact}, for every sequence $\{x_n \}_{n\in \mathbb{N}}\subseteq X$ such that $x_n\leq x_{n+1}$, one has that $\bigcap_{n\in \mathbb{N}} F_{x_n}\neq \emptyset$. Argue by contradiction assuming that $X$ is $\sigma$-compact. By Theorem \ref{thm_caracterizacion_sigma_compacto}, $X$ contains a countable cofinal subset $C=\{c_n \}_{n\in \mathbb{N}}$. We construct recursively an increasing sequence $\{x_n\}_{n\in\mathbb{N}}$ such that $x_n\geq c_n$ for every $n\in \mathbb{N}$. First, set $x_1=c_1$. Assume that $x_n$ is known. By Theorem \ref{conexion_caracterizacion}, $X$ is ultraconnected. Therefore, there exists $z_n$ such that $x_n,c_{n+1}\leq z_n$. Moreover, $X$ is non-compact and connected, consequently, $X$ does not have maximal points (see Proposition \ref{propiedadesbasicas}) and there exists $x_{n+1}$ such that $z_n<x_{n+1}$. Hence, $\{ x_n\}_{n\in \mathbb{N}}$ satisfies that $x_n<x_{n+1}$ and also that $c_n\leq x_n$ for every $n\in \mathbb{N}$. By hypothesis, $\bigcap_{n\in \mathbb{N}} F_{x_n}\neq \emptyset$. Let $t\in \bigcap_{n\in \mathbb{N}} F_{x_n}$. Then $t\geq  x_n$ for every $n\in \mathbb{N}$. On the other hand, $C$ is cofinal, which means that there exists $c_m\in C$ such that $c_m\geq t$. However,
    \[ c_m\leq x_m<x_{m+1}\leq t\]
    and hence $c_m< t$, contradicting $t\leq c_m$. Therefore, $X$ cannot be both $2$-pseudocompact and $\sigma$-compact.


\end{proof}

\begin{cor}\label{cor_excluyentes}
    Let $X$ be a non-compact Alexandroff paratopological group with finitely many connected components. Then $X$ cannot be both $2$-pseudocompact and $\sigma$-compact.
\end{cor}
\begin{proof}
    Suppose $X$ is $2$-pseudocompact. By Theorem \ref{thm_2-pseudo_finitos}, each connected component of $X$ is $2$-pseudocompact. By Theorem \ref{thm_excluyentes}, each connected component is not $\sigma$-compact. On the other hand, $X$ is $\sigma$-compact if and only if each connected component is $\sigma$-compact. This gives a contradiction if we assume that $X$ is also $\sigma$-compact.
\end{proof}

Although the preceding theorem shows that $2$-pseudocompactness and $\sigma$-compactness are mutually exclusive, the following example shows that they do not form a dichotomy.

\begin{ex}\label{ex_funciones_naturales_r}
    Let $\mathbb{R}^{\mathbb{N}}$ be the set of maps from the set of natural numbers $\mathbb{N}$ to the set of the real numbers. Declare that 
    \[ f\leq g \quad \textnormal{if and only if }\quad  f(n)\leq g(n) \quad \textnormal{for all }n\in \mathbb{N},\]
    where $f,g\in \mathbb{R}^{\mathbb{N}}$. This relation defines a partial order on $\mathbb{R}^{\mathbb{N}}$. It is evident that $\mathbb{R}^{\mathbb{N}}$ is an Abelian group with the operation given by $(f+g)(n)=f(n)+g(n)$. Moreover, this group operation is continuous and turns $\mathbb{R}^{\mathbb{N}}$ into a non-compact $T_0$ Alexandroff paratopological group. We verify that $\mathbb{R}^{\mathbb{N}}$ is hyperconnected. For every pair of points $f$ and $g$ we can consider $h$ defined as $h(n)=\max\{f(n),g(n) \}$, which satisfies that $f,g\leq h$, i.e., $X$ is a directed set. By Theorem \ref{conexion_caracterizacion}, we conclude that $X$ is hyperconnected. Suppose that there exists a countable cofinal subset $C=\{f_n \}_{n\in \mathbb{N}}$. Define $g\in \mathbb{R}^\mathbb{N}$ by 
    \[ g(n)=f_n(n)+1.\]
    Clearly, for each $n\in \mathbb{N}$, one has that $g(n)>f_n(n)$, which means that $g\nleq f_n$ for every $n\in \mathbb{N}$. Thus, $X$ is not $\sigma$-compact by Theorem \ref{thm_caracterizacion_sigma_compacto}.

    Consider now the sequence $\{f_n\}_{n\in \mathbb{N}}\subset \mathbb{R}^{\mathbb{N}}$ defined as 
    \[ f_n(m)=n \text{ for every } m\in \mathbb{N}. \]
    It is evident that $f_{n}<f_{n+1}$. Moreover, $\bigcap_{n\in \mathbb{N}}F_{f_n}=\emptyset$. If there exists $h\in \bigcap_{n\in \mathbb{N}}F_{f_n}$, then $h(m)\geq f_n(m)=n$ for every $n\in \mathbb{N}$, which is not possible. By Theorem \ref{thm_car_2_pseudocompact}, $X$ is not $2$-pseudocompact.
\end{ex}

\section{Classification by number of connected components}\label{sec_numero_componentes_conexas}

\subsection{$H$-compactness, pseudocompactness and star-compactness}

We characterize $H$-compactness, pseudocompactness and star-compactness within the setting of non-compact Alexandroff paratopological groups in terms of the number of connected components. 

\begin{df}
    A topological space $X$ is called {\it $H$-compact} if any open cover $\mathcal{U}$ of $X$ contains a finite subfamily $\mathcal{V}\subseteq\mathcal{U}$ such that $X =\bigcup_{V\in\mathcal{V}}\overline{V}$.
\end{df}

\begin{thm}\label{thm_car_h-compact}
    Let $X$ be a non-compact Alexandroff paratopological group. Then $X$ is $H$-compact if and only if $X$ has finitely many connected components.
\end{thm}
\begin{proof}
    Suppose that $X$ has a finite number of connected component: $X_1,\ldots,X_n$. Given an open cover $\mathcal{U}$ of $X$, consider for each $X_i$ an element $U_i\in \mathcal{U}$ such that $X_i\cap U_i \neq\emptyset$. Let $\mathcal{V}=\{U_i \}_{i=1,\ldots,n}$. By Theorem \ref{dense}, $X_i\cap U_i$ is dense in $X_i$. Thus, $X=\bigcup_{i=1}^n \overline{U_i}$.

    Suppose that $X$ is H-compact. Consider the open cover $\mathcal{U}=\{ U_x\}_{x\in X}$. Then there exists a finite set $M\subset X$ such that $X=\bigcup_{x\in M} \overline{U_x} $. It is evident that each $x\in M$  belongs to just one connected component. Moreover, by Theorem \ref{conexion_caracterizacion}, $\overline{U_x}$ provides one connected component. Since $X=\bigcup_{x\in M}\overline{U_x}$, there are at most $|M|$ connected components.
\end{proof}

\begin{df}
    A topological space $X$ is said to be {\it pseudocompact} if every real continuous function defined on it is bounded.
\end{df}

\begin{thm}\label{car_pseudocompact}
    Let $X$ be a non-compact Alexandroff paratopological group. Then $X$ is pseudocompact if and only if $X$ has finitely many connected components.
\end{thm}
\begin{proof}
    By Proposition \ref{propiedadesbasicas}, every continuous map $f:X\rightarrow \mathbb{R}$ is constant on the connected components of 
    $X$. From this, we deduce the desired result.
\end{proof}

\begin{df}
Let $X$ be a topological space and let $\mathcal{U}$ be a family of subsets of $X$. For every $A\subseteq X$, the \emph{star of $A$ with respect to $\mathcal{U}$} is defined as \[ \operatorname{St}(A,\mathcal{U}) = \bigcup\{U\in\mathcal{U}:U\cap A\neq\emptyset\}. \] We say that $X$ is \emph{star-compact} if, for every open cover $\mathcal{U}$ of $X$, there exists a finite subset $F\subseteq X$ such that \[ \operatorname{St}(F,\mathcal{U})=X. \]
\end{df}

\begin{thm}
Let $X$ be a connected non-compact Alexandroff paratopological group and let $x\in X$. If $\mathcal{U}$ is an open cover, then $\operatorname{St}(\{ x\},\mathcal{U})=X.$ Particularly, $X$ is star-compact.
\end{thm}
\begin{proof}
    By Theorem \ref{conexion_caracterizacion}, $X$ is hyperconnected. Therefore, for every $y\in X$, there exists $z_y$ such that $z_y\geq x,y$. Moreover, $\mathcal{U}$ is an open cover, so there exists $U\in \mathcal{U}$ such that $z_y\in U$. Consequently, $x\in U$ because $U$ is a lower set. From here, we deduce the desired result.
\end{proof}

\begin{cor}\label{char_star_compact}
    Let $X$ be a non-compact Alexandroff paratopological group. Then $X$ is star-compact if and only if $X$ has a finite number of connected components.
\end{cor}

By combining Theorem \ref{thm_car_h-compact}, Theorem \ref{car_pseudocompact} and Corollary \ref{char_star_compact} we have the following characterization.
\begin{thm}
    Let $X$ be a non-compact Alexandroff paratopological group. Then the following conditions are equivalent:
    \begin{enumerate}
        \item $X$ is pseudocompact.
        \item $X$ is $H$-compact.
        \item $X$ is star-compact.
        \item $X$ has a finite number of connected components.
    \end{enumerate}
\end{thm}

\subsection{Limit compactness, weakly Lindelöf and countable chain condition}

As we did in other subsections, we start by recalling the definitions that are needed and then provide characterizations. 

\begin{df}
    A topological space $X$ is said to be \emph{limit compact}  if every infinite subset of $X$ has a limit point in $X$. A point $p$ is a \emph{limit point} of a set $A$ if every open set containing $p$ contains at least one point of $A$ distinct from $p$. 
\end{df}
\begin{thm}\label{limit_compact}
    Let $X$ be a non-compact Alexandroff paratopological group. Then $X$ is limit compact if and only if $\mathrm{width}(X)<\infty$.
\end{thm}
\begin{proof}
    Suppose that $X$ is limit compact. We argue by contradiction assuming that $\mathrm{width}(X)=\infty$. Then, we can consider an antichain $A$ such that $|A|=\infty$. For every $a\in A$, one has that $U_a\cap A=\{ a\}$, otherwise we obtain a contradiction with the fact that $A$ is an antichain. Thus, $A$ does not have a limit point in $X$.

    Suppose that $\mathrm{width}(X)<\infty$. This gives that there are no antichains with infinite points. Consquently, if $A$ is an infinite subset, then there exist $a,b\in A$ such that $a<b$. This gives that $b$ is a limit point. 
\end{proof}
\begin{df}
    A topological space $X$ satisfies the \emph{countable chain condition} if every family of pairwise disjoint open subsets of $X$ is countable.
\end{df}
\begin{thm}\label{chain_condition}
    Let $X$ be a non-compact Alexandroff paratopological group. Then $X$ satisfies the countable chain condition if and only if $X$ has a countable number of connected components.
\end{thm}
\begin{proof}
    By Theorem \ref{conexion_caracterizacion}, each connected component is a hyperconnected topological space. If $\mathcal{U}$ is a family of pairwise disjoint open subsets, then each member of $\mathcal{U}$ belongs only to one of the connected components of $X$. From here, we may deduce the desired result.
\end{proof}

\begin{df}
     A topological space $X$ is \emph{weakly Lindelöf} if its weak Lindelöf number is $\aleph_0$, that is, if every open covering of $X$ contains a countable subfamily whose union is dense in $X$. 
\end{df}

\begin{thm}\label{weakly_lindelof}
    Let $X$ be a non-compact Alexandroff paratopological group. Then $X$ is weakly Lindelöf if and only if $X$ has a countable number of connected components.
\end{thm}
\begin{proof}
    Every open subset of a connected non-compact Alexandroff paratopological group is dense (see Theorem \ref{conexion_caracterizacion}). From this, it is straighforward to deduce the result.
\end{proof}

Combining Theorem \ref{limit_compact}, Theorem \ref{chain_condition} and Theorem \ref{weakly_lindelof}, we deduce the following results:

\begin{thm}
    Let $X$ be a non-compact Alexandroff paratopological group. If $X$ is limit compact, then $X$ is weakly Lindelöf and satisfies the countable chain condition.
\end{thm}
\begin{proof}
    We argue by contradiction. Suppose that $X$ is not wealky Lindelöf or does not satisfy the countable chain condition. Then $X$ has uncountable connected components $X_i$ with $i\in I$ (see Theorem \ref{chain_condition} and Theorem \ref{weakly_lindelof}). For each connected component choose $x_i\in X_i$ and set $A=\{x_i \}_{i\in I}$. The set $A$ is clearly an infinite antichain, therefore, $\mathrm{width}(X)=\infty$, which entails a contradiction with Theorem \ref{limit_compact}.
    \end{proof}
\begin{thm}
    Let $X$ be a non-compact Alexandroff paratopological group. Then the following conditions are equivalent:
    \begin{enumerate}
        \item $X$ satisfies the countable chain condition.
        \item $X$ is weakly Lindelöf.
        \item $X$ has a countable number of connected components.
    \end{enumerate}
\end{thm}

\section{Properties that do not hold}\label{sec_fallan}

In this subsection, we prove properties that do not hold for non-compact Alexandroff paratopological groups. We start by proving a useful lemma in what follows.

\begin{lem}\label{lem_clausura_es_no_compacta}
    Let $X$ be a connected non-compact Alexandroff space and $A\subseteq X$, then the closure of $A$ is not compact.
\end{lem}
\begin{proof}
Let $e$ denote the identity element of $X$. By hypothesis, there exists $y\in X$ such that $e< y$. By the continuity of the multiplication map, we obtain the following chain:
\[e<y<y^2<\cdots <y^n<y^{n+1}\cdots \]
On the other hand, assume by contradiction that $F_A$, i.e., the closure of $A$, is compact. Then there exists a finite subset $M\subseteq F_A$ such that for every $y\in F_A$, there exists $m_y\in M $ with $y<m_y$ (see Theorem \ref{thm_caracterizacion_compactos}).
By multiplying the previous chain by $m$ for some $m\in M$, we obtain:
\[ m<my<my^2<\cdots <my^n<my^{n+1}\cdots\]
Clearly, $my^n\in F_A$ for every $n\in \mathbb{N}$ because the closure of $A$ is an upper set. This entails a contradiction with the finiteness of $M$.
\end{proof}

\begin{df}
    A topological space $X$ is {\it $\omega$-bounded} if each countable subset of $X$ has compact closure.
\end{df}
\begin{df}
    A topological space $X$ is {\it totally countably compact} if each sequence of $X$ contains a subsequence with compact closure.
\end{df}
\begin{df}
    A topological space $X$ is called {\it strongly locally compact} if each point can be contained in an open set whose closure is compact. 
\end{df}


\begin{thm}\label{thm_no_funciona}
    Let $X$ be a  non-compact Alexandroff paratopological group. Then
    \begin{enumerate}
    \item $X$ is not $\omega$-bounded.
    \item $X$ is not totally countably compact.
    \item $X$ is not strongly locally compact if there $X$ does not have singleton connected components.
    \end{enumerate}
\end{thm}
\begin{proof}

1.-3. If $X$ is connected, then it is straightforward from Lemma \ref{lem_clausura_es_no_compacta}.  If $X$ is not connected, then we may apply Lemma \ref{lem_clausura_es_no_compacta} to each connected component.

\end{proof}

\begin{rem}
    Note that by removing the extra hypothesis from 3. Theorem \ref{thm_no_funciona}, the result does not hold. For instance, consider the group of integer numbers $\mathbb{Z}$ with the discrete topology. Then $\mathbb{Z}$ is trivially strongly locally compact.
\end{rem}

\begin{df}
    A toplogical space $X$ is {\it sequentially feebly compact} if for every sequence $\{U_n\}_{n\in \mathbb{N}}$ of nonempty open sets of $X$ there exists
a point $x \in X$ and an infinite set $I \subseteq \mathbb{N}$ such that for every neighborhood $O_x \subseteq X$ of $x$ the set
$\{n \in I : O_x \cap U_n \neq \emptyset \}$ is finite.

\end{df}

\begin{prop}\label{no_secuencial_feebly_compact} 
    Let $X$ be a  non-compact Alexandroff paratopological group. Then 
    $X$ is not sequentially feebly compact.
\end{prop}

\begin{proof}
Suppose that $X$ is connected and let $\{U_n\}_{n\in\mathbb{N}}$ be a sequence of nonempty open subsets of $X$. By Theorem \ref{every_paratopological_alexandroff_hyperconexosiiiiii}, every neighborhood $O_x$ of $x\in X$ intersects each $U_n$. That is, \[ O_x\cap U_n\neq\emptyset \] for all $n\in\mathbb{N}$. If $X$ is not connected, the same argument applies to each connected component of $X$.
\end{proof}

\begin{df}
    A paratopological group $X$ is {\it precompact} if  for each neighborhood $U$ of the identity of $X$, there exists a finite subset $F$ of $X$ such that $FU = X$.
\end{df}

\begin{prop}
    Let $X$ be a non-compact Alexandroff paratopological group. Then $X$ is not precompact.
\end{prop}

\begin{proof}
Suppose that $X$ is precompact. Then there exists a finite subset $F$ of $X$ such that $FU_e=X$ where $e$ denotes the identity element. However, by the continuity of the multiplication map, for every $x\in FU_e$ there exists $y\in F$ such that $y\geq x$. Then $X=\bigcup_{x\in F} U_x$. This gives that $X$ is compact because $X$ would have a finite number of maximal points (see Theorem \ref{thm_caracterizacion_compactos}), which entails a contradiction.
\end{proof}

\begin{df}
    A topological space is called {\it (countable) metacompact}  if every (countable) open cover has an open point finite refinement. A cover is point finite if each point belongs to only finitely many sets in the covering.
\end{df}

\begin{prop}
    Let $X$ be a non-compact Alexandroff paratopological group with no singleton connected components. Then $X$ is not  (countable) metacompact.
\end{prop}
\begin{proof}
Let $e$ denote the identity element and assume $X$ is connected. For every $x\in X$ there exists $d_x\in X$ such that $d_x\geq x,e$ because $X$ is connected (see Theorem \ref{conexion_caracterizacion}). Consider the open cover $\mathcal{U}=\{U_{d_x} \}_{x\in X}$. Suppose by contradiction that $X$ is metacompact. Therefore, there exists an open  point finite refinement $\mathcal{V}$. Consequently, the identity element $e$ belongs to only finitely many sets in the covering. However, every open set of $\mathcal{U}$ contains the identity element $e$. This implies that $\mathcal{V}$ is finite, which leads to a contradiction with the non-compactness of $X$ (see Theorem \ref{thm_caracterizacion_compactos}). If $X$ is not connected, we may use the same arguments to the connected component that contains the identity element.
\end{proof}

\section{Covering dimensions}\label{sec_covering_dimension}

In this section, we recall the classical notion of covering dimension and several of its generalizations, and study them in the context of the paratopological groups considered in this paper. We also highlight some connections with the notions introduced in the preceding sections.
\begin{df}
   Let $X$ be a topological space. An open cover $\mathcal{U}$ of $X$ is said to have \emph{order at most $n$} if every point of $X$ belongs to at most $n+1$ members of $\mathcal{U}$, or equivalently, if the intersection of any $n+2$ pairwise distinct members of $\mathcal{U}$ is empty. We say that $X$ has \emph{Lebesgue covering dimension at most $n$}, and write $\dim(X)\leq n$, if every open cover of $X$ admits an open refinement of order at most $n$. The \emph{Lebesgue covering dimension} of $X$ is the least nonnegative integer $n$ for which this property holds. If no such integer exists, we write $\dim(X)=\infty$ and say that $X$ is \emph{infinite-dimensional}.

\end{df}
First, we establish a general result within the theory of Alexandroff spaces. 

\begin{thm}\label{covering_dimension_implica_maximales}
    Let $X$ be an Alexandroff space. If $\dim (X)<\infty$, then $X$ contains a subset of maximal points $M$ such that for every $x\in X$ there exists $m\in M$ with $x\leq m$.
\end{thm}
\begin{proof}
    Suppose $\dim(X)=n$ and argue by contradiction by assuming that there are no maximal points. Consider the open cover $\mathcal{U}=\{U_x \}_{x\in X}$. Then $\mathcal{U}$ admits an open refinement $\mathcal{V}$ such that for every $x\in X$, $x$ is contained in at most $n+1$ members of $\mathcal{V}$, i.e., there exists $M\subset X$ such that $x\leq y$ for every  $y\in M$ and $|M|\leq n+1$. For each $y\in M$, we may find $z\in X$ such that $y<z$ because there are no maximal points in $X$. Since $\mathcal{V}$ is also an open cover, $U_z\in \mathcal{V}$ for some $z>y$. On the other hand, $x\in U_z$. We may repeat this argument until we reach a contradiction with the fact that $\dim(X)=n$.
\end{proof}

The purpose of the following examples is twofold. On the one hand, we prove that $\dim(X)=n$ does not imply that there are $n$ maximal points. On the other hand, we show that the coverse of Theorem \ref{covering_dimension_implica_maximales} does not hold.
\begin{ex}
     Let $X=\{ x_i\}_{i\in \mathbb{N}}$ such that $x_i<x_{i+1}$ if $i$ is odd and $x_{i}>x_{i+1}$ if $i$ is even. It is straightforward to verify that $\dim(X)=1$ and that there is a countable number of maximal points.
    
\end{ex}
\begin{ex}
    Let $X=\mathbb{N}\cup \{ *\}$ with $*<n$ for every $n\in \mathbb{N}$. Then $M=\mathbb{N}$ is a set of maximal points such that for every $x\in X$ there exists $m\in M$ with $x\leq m$. However, $\dim (X)=\infty$.
\end{ex}

\begin{cor}
    Let $X$ be a non-compact Alexandroff paratopological group. Then $X$ is infinite-dimensional.
\end{cor}
\begin{proof}
It is straightforward from Theorem \ref{covering_dimension_implica_maximales} and the fact that there are no maximal points in $X$ (see Proposition \ref{propiedadesbasicas}).
\end{proof}

\begin{df}
    A topological space $X$ is said to have \emph{local covering dimension at most $n$} if every point $x\in X$ has an open neighborhood $U$ such that $ \dim (U)\leq n$.
\end{df}

\begin{prop}\label{prop_localmente_finit_covering_dimension}
    Let $X$ be a non-compact Alexandroff paratopological group. Then $X$ has local covering dimension $0$. Particularly, for every $x\in X$, $\dim(U_x)=0$.
\end{prop}
\begin{proof}
    Let $x\in X$. Then $U_x$ is an element in every open cover $\mathcal{U}_x$ of $U_x$. Hence, $\mathcal{V}=\{ U_x\}$ is an open refinement of  $\mathcal{U}_x$  with order at most $0$. 
\end{proof}

\begin{rem}
    Non-compact Alexandroff paratopological groups have infinite Lebesgue covering dimension, although its local covering dimension is zero at every point.
\end{rem}

\begin{df}
    A topological space $X$ is said to be \emph{countable-dimensional} if there exists a countable family $\{X_n\}_{n\in\mathbb{N}}$ of finite dimensional subspaces of $X$ such that \[ X=\bigcup_{n\in\mathbb{N}}X_n. \]
\end{df}

\begin{thm}\label{sigma_compacto_implica_countable}
    Let $X$ be a non-compact Alexandroff paratopological group. If $X$ is $\sigma$-compact, then $X$ is countable-dimensional. 
\end{thm}
\begin{proof}
    If $X$ is $\sigma$-compact, then there exists a countable cofinal subset $M$ in $X$ (see Theorem \ref{thm_caracterizacion_sigma_compacto}). By Proposition \ref{prop_localmente_finit_covering_dimension}, $\dim(U_x)=0$ for every $x\in M$. Hence, $\{U_x\}_{x\in M}$ is a countable family of finite dimensional subspaces of $X$ satisfying $X=\bigcup_{x\in M}U_x$.   
\end{proof}

The following example demonstrates that the opposite of Theorem \ref{sigma_compacto_implica_countable} does not hold.

\begin{ex}\label{ex_ZxR}
    Let us consider the non-compact Alexandroff paratopological group of the integer numbers $\mathbb{Z}$ (see Example \ref{ex_z_sigma_compacto}). Let $\mathbb{R}$ be the group of the real numbers with the discrete topology. Consider the non-compact Alexandroff paratopological group $\mathbb{Z}\times \mathbb{R}$ obtained by the product of $\mathbb{Z}$ and $\mathbb{R}$. It is easy to verify that $\mathbb{Z}\times \mathbb{R}$ is not $\sigma$-compact because there is no countable cofinal subset (see Theorem \ref{char_sigma_compact}). Let $n\in\mathbb{N}$, consider $X_n$ defined by
    \[
    X_n=\bigcup_{a\in \mathbb{R}} U_{(n,a)}.
    \]
    By Proposition \ref{prop_localmente_finit_covering_dimension}, $\dim(U_{(n,a)})=0$. Moreoever, $X=\bigcup_{n\in \mathbb{N}} X_n$, which yields that $\mathbb{Z}\times \mathbb{R}$ is countable-dimensional.
    \end{ex}

\begin{df}
    A topological space $X$ is \emph{strongly countable-dimensional} if $ X=\bigcup_{n\in\mathbb{N}} F_n$, where each $F_n$ is a closed finite-dimensional subspace of $X$.
\end{df}

\begin{prop}
    Let $X$ be a non-compact Alexandroff paratopological group. Then $X$ is not strongly countable-dimensional.
\end{prop}
\begin{proof}
    A closed set in $X$ does not have maximal points and, therefore, $X$ is not strongly countable-dimensional (see Theorem \ref{covering_dimension_implica_maximales}) . 
\end{proof}

\section{Consequences and relations}\label{sec_consequencias}

We devote this section to obtain consequences and implications between the notions described in the previous sections.

\begin{thm}
Let $X$ be a non-compact Alexandroff paratopological group with infinitely many connected components. Then $X$ is neither $2$-pseudocompact nor countably $o$-compact, and it is not a Baire space.
\end{thm}
\begin{proof}
    Let $\{ X_n\}_{n\in \mathbb{N}}$  be a set of connected components of $X$. Consider $U_n=\bigcup_{i\geq n} X_i$. Then $\{U_n \}_{n\in \mathbb{N}}$ is a decreasing sequence of nonempty open sets. However, $\bigcap_{n\in \mathbb{N}} U_n$ is the empty set. This shows that $X$ is not countably $o$-compact. We can consider the same sequence to prove that $X$ is neither $2$-pseudocompact nor a Baire space.
\end{proof}

\begin{thm}
    Let $X$ be a non-compact Alexandroff paratopological group. If $X$ has an uncountable number of connected components, then
    \begin{itemize}
    \item[1.] $X$ is not $\sigma$-compact, $\sigma$-locally compact neither Lindelöf
    
    \item[2.] $X$ is not pseudocompact, $H$-compact neither star-compact.
    \item[3.] $X$ is not weakly Lindelöf and does not satisfy the countable chain condition.
    \end{itemize}
\end{thm}

\begin{thm}
    Let $X$ be a non-compact Alexandroff paratopological group. If $X$ is not countable-dimensional, then
    \begin{itemize}
    \item[1.] $X$ is not $\sigma$-compact, $\sigma$-locally compact neither Lindelöf
    
    \item[2.] $X$ has an uncountable cofinal subset.
    
    \end{itemize}
\end{thm}

\begin{thm}\label{sigmacompactoimplicaweaklylindelof}
    Let $X$ be a non-compact Alexandroff paratopological group. If $X$ is $\sigma$-compact then it is weakly Lindelöf.
\end{thm}
\begin{proof}
    By Theorem \ref{char_sigma_compact}, there exists a countable cofinal subset and it is straightforward to prove that this implies $X$ has a countable number of connected component. By Theorem \ref{weakly_lindelof}, we conclude the proof.
\end{proof}

\begin{ex}\label{ex_weakly_lindelof_no_implica_sigma_compact}
    Let $\mathbb{Z}^\mathbb{N}$ the paratopological group from Example \ref{ex_2-pseudocompact} and let $\mathbb{Z}$ be the group of integer numbers with the discrete topology. Then $\mathbb{Z}^\mathbb{N}\times \mathbb{Z}$ is weakly Lindelöf because it has a countable number of connected component. However, it is immediate to verify that $\mathbb{Z}^\mathbb{N}\times \mathbb{Z}$ is not $\sigma$-compact.
\end{ex}
\begin{thm}\label{limit_compact_sigma_compact_chain}
    Let $X$ be a non-compact Alexandroff paratopological group. If $X$ is limit compact and every chain in $X$ is $\sigma$-compact, then $X$ is $\sigma$-compact.
\end{thm}
\begin{proof}
    By hypothesis and Theorem \ref{limit_compact}, we have $\mathrm{width}(X)=n<\infty$. By Dilworth's theorem, $X$ can be partitioned into $n$ chains. Moreover, each of these chains is $\sigma$-compact and therefore admits a countable cofinal subset. Let $M$ be the union of these countable cofinal subsets. Since $M$ is a finite union of countable sets, it is countable. Furthermore, $M$ is clearly cofinal in $X$. Hence, $X$ is $\sigma$-compact (see Theorem \ref{char_sigma_compact}).
\end{proof}

With the following example we demonstrate that we cannot remove from Theorem \ref{limit_compact_sigma_compact_chain} the hypothesis of limit compactness.
\begin{ex}\label{ex_hipotesis_limit_compact_importa}
    Let us consider the paratopological group $\mathbb{Z}\times \mathbb{R}$ from Example \ref{ex_ZxR}.  Clearly, every chain in  $\mathbb{Z}\times \mathbb{R}$ is $\sigma$-compact. However, the width of  $\mathbb{Z}\times \mathbb{R}$ is uncountable, then by Theorem \ref{limit_compact}, $\mathbb{Z}\times \mathbb{R}$ is not limit compact. Moreover, this group is not $\sigma$-compact.
\end{ex}
\begin{rem}
    In the proof of Theorem~\ref{limit_compact_sigma_compact_chain}, we use Dilworth's theorem. It may seem natural to seek a similar result by replacing limit compactness, and hence the finiteness of the width, with the assumption that the space has countably infinite width. However, Dilworth's theorem does not, in general, extend in this direction: a poset of countably infinite width need not admit a partition into countably many chains.
\end{rem}

\bibliography{bibliografia}
\bibliographystyle{plain}

\end{document}